\documentclass[preprint,12pt]{elsarticle}

\newtheorem{theorem}{Theorem}[section]
\newtheorem{corollary}{Corollary}[section]

\usepackage{amssymb}
\begin{document}

\begin{frontmatter}

\title{ Berry-Esseen bounds for self-normalized martingales}
\author[cor1]{Xiequan Fan}
%\author[cor2]{Ion Grama}
%\author[cor2]{Quansheng Liu}
\author[cor3]{Qi-Man Shao }
\address[cor1]{Center for Applied Mathematics,
Tianjin University, Tianjin 300072,  China}
%\address[cor2]{Universit\'{e} de Bretagne-Sud, LMBA, UMR CNRS 6205,
% Campus de Tohannic, 56017 Vannes, France}
\address[cor3]{Department of Statistics, The Chinese University of Hong Kong, Shatin, NT, Hong Kong}

\begin{abstract}
A Berry-Esseen bound is obtained for self-normalized martingales under the assumption of  finite moments. The   bound
coincides with the classical Berry-Esseen bound for standardized martingales. An example is given to show the optimality of the bound.  Applications to Student's  statistic  and  autoregressive process are also discussed.
\end{abstract}

\begin{keyword} Self-normalized process, Berry-Esseen bounds, martingales, Student's  statistic,  autoregressive process
\vspace{0.3cm}
\MSC Primary 60G42; 60F05;  Secondary 60E15
\end{keyword}

\end{frontmatter}

% \linenumbers

\section{Introduction}
Let $(X_i)_{i\geq 1}$ be a sequence of independent non-degenerate real-valued random variables  with zero means,
and let $$S_n=\sum_{i=1}^n X_i \ \ \    \textrm{and} \ \ \     V_n^2=\sum_{i=1}^n X_i^2$$
be the partial sum and the partial quadratic sum,  respectively.
The self-normalized sum is defined as $S_n/V_n.$  The study of the asymptotic behavior of self-normalized sums
has a long history. When $(X_i)_{i\geq 1}$ are i.i.d. in the domain of normal and stable law,  Logan et al.\ \cite{LMRS73} obtained the weak convergence for the self-normalized sum, while
Gin\'{e} et al.\  \cite{GGM97} proved that $S_n/V_n$ is asymptotically normal if and only if $X_1$ belongs to
the domain of attraction of a normal law.
 Under the same  necessary and sufficient condition, Cs\"{o}rg\H{o} et al.\ \cite{CSW03} proved
a self-normalized type Donsker's theorem.  For general independent random variables with finite $(2+\delta)^{th}$ moments, where $0 < \delta \leq 1$,
 Bentkus, Bloznelis and G\"{o}tze \cite{BBG96} (see also  Bentkus and   G\"{o}tze \cite{BG96} for i.i.d.\,case)
 have obtained the following Berry-Esseen bound : If $\mathbf{E} |X_i|^{2+\delta}    < \infty$ for $\delta\in (0, 1],$ then these exists an absolute constant $C$ such that
\begin{eqnarray*}
\sup_x \Big|\mathbf{P}(S_n/ V_n \leq x) - \Phi(x) \Big|  &\leq& C\, B_n^{-2-\delta}\sum_{i=1}^n\mathbf{E}|X_i|^{2+\delta}   ,
\end{eqnarray*}
 where $B_n^2 = \sum_{i=1}^n \mathbf{E}X_i^2$, and $\Phi \left( x\right)$ is the standard normal distribution function.
%\begin{eqnarray*}
%&&\sup_x \Big|\mathbf{P}(S_n/ V_n \leq x) - \Phi(x) \Big|  \nonumber   \\
%&&\ \ \ \ \ \ \ \ \ \ \ \  \leq \, C\, \Big( B_n^{-2}\sum_{i=1}^n\mathbf{E}[X_i^2\mathbf{1}_{\{|X_i|>B_n\}} ] + B_n^{-3}\sum_{i=1}^n\mathbf{E}[|X_i|^3\mathbf{1}_{\{|X_i|\leq B_n\}} ]  \Big),
%\end{eqnarray*}
%In particular, the last inequality implies that
It is worth noting that the last bound coincides with the classical Berry-Esseen bound for  standardized  partial sums $S_n/ B_n$ and it is the best possible.
For the related error  of $\mathbf{P}(S_n/ V_n \geq x)$ to $1-\Phi(x),$ we refer to
 Shao \cite{S99}, Jing, Shao and  Wang \cite{JSW03}. In these papers, self-normalized Cram\'{e}r type moderate deviation theorems have been established. We also refer to de la Pe\~{n}a, Lai and Shao \cite{PLS09}, Shao and Wang \cite{SW13} and Shao and Zhou \cite{SZ17} for surveys on recent developments on self-normalized limit theory.

 Despite the fact that the case  for self-normalized sums of independent random variables  is well studied, we are not aware of Berry-Esseen bounds for self-normalized martingales in the literature. The main purpose of this paper is to fill this gap.

We first recall some Berry-Esseen bounds for standardized martingale difference sequence.
Let $(X _i,\mathcal{F}_i)_{i=0,...,n}$ be a finite  sequence of martingale differences defined on a
 probability space $(\Omega ,\mathcal{F},\mathbf{P})$,  where $X
_0=0 $ and $\{\emptyset, \Omega\}=\mathcal{F}_0\subseteq ...\subseteq \mathcal{F}_n\subseteq
\mathcal{F}$ are increasing $\sigma$-fields.   Set
\begin{equation}\label{matingal}
S_{0}=0,\ \ \ \ \ \ \ S_k=\sum_{i=1}^kX _i,\quad k=1,...,n.
\end{equation}
Then $S=(S_k,\mathcal{F}_k)_{k=0,...,n}$ is a martingale.
Let $[ S]$ and $\langle S \rangle$  be, respectively, the squared variance and the conditional variance  of the
martingale $S,$ that is
\begin{eqnarray*}
[S]_0=0,\ \ \ \ \ [ S]_k=\sum_{i=1}^k X _i^2
\end{eqnarray*}
and
\begin{eqnarray*}
\langle S\rangle_0=0,\ \ \ \ \ \langle S \rangle_k=\sum_{i=1}^k \mathbf{E} [ X _i^2  |  \mathcal{F}_{i-1} ]  ,\quad k=1,...,n.
\end{eqnarray*}
Suppose that   $ \mathbf{E} |X_i|^{2p}    < \infty $ for some $p>1 $ and all $i=1,...,n.$
Define
\begin{equation}\label{defNn}
N_n = \sum_{i=1}^n \mathbf{E} |X_i|^{2p}  + \mathbf{E} \big| \langle S\rangle_n-1\big|^p.
\end{equation}
When $p\in (1, 2],$ Heyde and Brown \cite{HB70} (see also Theorem 3.10 of Hall and Heyde \cite{HH80}) proved that  there exists a constant $C_p$  depending only on $p$   such that
\begin{equation}\label{sfdf}
 \sup_{ x \in \mathbf{R}}\Big|\mathbf{P}(S_n  \leq x)-\Phi \left( x\right) \Big|  \leq  C_p \, N_n^{1/(2p+1) }.
\end{equation}
Later, Haeusler \cite{H88}
gave an extension of (\ref{sfdf}) to all $p\in (1, \infty)$. Moreover, Haeusler also gave an example to justify that his bound
is asymptotically the best possible. It is remarked that the  $(X_i)_{1 \leq i \leq n}$ is standardized, that is,
$\sum_{i=1}^n \mathbf{E}X_i^2$ is close to $1$.

In this paper, we prove that the Berry-Esseen bound  (\ref{sfdf})  also holds for self-normalized martingales $S_n/ \sqrt{[S]_n}$ and normalized martingales $S_n/ \sqrt{\langle S\rangle_n}$.  Moreover, we also justify the optimality of our bounds. Applications to Student's  statistic  and  autoregressive process  are   discussed.

The paper is organized as follows. Our main results are stated and discussed in Section \ref{sec2}.
The applications   are given in Section  \ref{secaps}.
Proofs of theorems are deferred to Section \ref{sec3}.

\section{Main results}\label{sec2}
The following theorem gives a counterpart of Haeusler's result \cite{H88}  for self-normalized martingales.
\begin{theorem}\label{th1}
Suppose that   $ \mathbf{E} |X_i|^{2p}    < \infty $ for some $p>1 $ and all $i=1,...,n.$
 Then there exists a constant $C_p$    depending only  on  $p  $  such that
\begin{equation}\label{ineq10}
\sup_{ x \in \mathbf{R}}\Big|\mathbf{P}\Big( \frac{S_n}{\sqrt{[S]_n}}  \leq x \Big)-\Phi \left( x\right) \Big|  \leq  C_p \, N_n^{1/(2p+1) },
\end{equation}
where $N_n$  is defined by (\ref{defNn}). Moreover, there exist a sequence of martingale differences $(X _i,\mathcal{F}_i)_{i=0,...,n} $ and a
positive constant $c_p$  depending only  on  $p$  such that
\begin{equation}\label{ineq07}
\sup_{ x \in \mathbf{R}}\Big|\mathbf{P}\Big( \frac{S_n}{\sqrt{[S]_n}}  \leq x \Big)-\Phi \left( x\right) \Big|  \, N_n^{-1/(2p+1) } \geq c_p.
\end{equation}
\end{theorem}

Clearly, inequality (\ref{ineq07}) shows that the bound (\ref{ineq10})
is asymptotically the best possible.

For a stationary martingale difference sequence, the term  $  \sum_{i=1}^n \mathbf{E}  |X_i|^{2p}   $ is of order $n^{ 1-p }$.
 Then  inequality (\ref{ineq10}) implies the following corollary.
\begin{corollary}\label{co01} Let $(X _i,\mathcal{F}_i)_{i \geq 1} $ be a stationary martingale difference sequence. Suppose that $  \mathbf{E}|X_1|^{2p} < \infty $ for some $p>1.$
Then there exists a constant $c_p, $  which does not depend on $n$,   such that
\begin{equation}\label{co01}
  \sup_{ x \in \mathbf{R}}\Big|\mathbf{P}\Big( \frac{S_n}{\sqrt{[S]_n}}  \leq x \Big)-\Phi \left( x\right) \Big|  \leq  c_p \, \Big(n^{1- p}+ \mathbf{E}  \big| \langle S\rangle_n-1\big|^p   \Big )^{1/(2p+1) }  .
\end{equation}
\end{corollary}

The next theorem gives a Berry-Esseen bound for normalized martingales  $S_n/ \sqrt{\langle S\rangle_n}$.

\begin{theorem}\label{th4}
Under the assumptions of Theorem \ref{th1},  the inequalities (\ref{ineq10}) and (\ref{ineq07}) hold when
 $S_n/\sqrt{[S]_n}   $ is replaced by $S_n/ \sqrt{\langle S \rangle_n} .$
\end{theorem}

For a stationary martingale difference sequence, the following result is a consequence of the last theorem.
\begin{corollary}
Assume the conditions of Corollary \ref{th1}. Inequality (\ref{co01})   holds when
 $S_n/\sqrt{[S]_n}   $ is replaced by $S_n/ \sqrt{\langle S \rangle_n} .$
\end{corollary}

\section{Applications} \label{secaps}
\subsection{Application to Student's $t$-statistic}
The study of self-normalized  partial sums
originates from Student's $t$-statistic. The Student's $t$-statistic $T_n$ is defined by
\[
T_n = \sqrt{n} \, \overline{X}_n / \widehat{\sigma},
\]
where $$\overline{X}_n = \frac{S_n}{n}  \ \ \ \textrm{and}\  \ \ \widehat{\sigma}^2 = \sum_{i=1}^n  \frac{(X_i - \overline{X}_n )^2}{ n-1}  .$$
It is known that for all $x\geq0,$
\[
\mathbf{P}\Big( T_n  > x \Big) = \mathbf{P}\Bigg(  \frac{S_n }{\sqrt{[S]_n}}  > x \Big(\frac{n}{n+x^2-1} \Big)^{1/2}  \Bigg ).
\]
 When $(X_i)_{i\geq1}$ is a sequence of  i.i.d.\ random variables,  Bentkus and   G\"{o}tze \cite{BG96} proved that
 if $\mathbb{E}|X_{i}| ^{2+\delta} < \infty$ for all $i=1,...,n$ and some $\delta \in (0, 1],$ then
\begin{eqnarray}
\sup_{x \in \mathbf{R}} \Big|\mathbf{P}(T_n \leq x) - \Phi(x) \Big| %% &\leq& C\, B_n^{-p}\sum_{i=1}^n\mathbf{E} |X_i|^p    \nonumber\\
&=& O\Big(n^{-\delta/2}\Big).
\end{eqnarray}

For martingales, we have the following analogue.
\begin{corollary} Let $(X _i,\mathcal{F}_i)_{i\geq 1} $ be a stationary martingale difference sequence.
Suppose that $  \mathbf{E}|X_1|^{2p}   < \infty $ for some $p>1.$
Then there exists a constant $C_p, $ which does not depend on $n$,   such that (\ref{co01})
holds  when $\mathbf{P}(S_n/\sqrt{[S]_n} \leq x)  $ is replaced by $\mathbf{P}(T_n \leq  x)  .$
\end{corollary}

%\subsection{Weak invariance principles}
%
%
%To be added

\subsection{Application to  autoregressive process}
Consider the autoregressive process given by
$$Y_{n+1} =  \theta Y_n + \varepsilon_{n+1}, \  \  n \geq 0, $$
where $Y_n$ and $\varepsilon_n$ represent  the observation and the driven noise, respectively.
The parameter $\theta $ is unknown and needs to be estimated at stage $n$ from the data
$Y_i, i\leq n.$ For sake of simplicity, we assume that $Y_0=0.$
We also assume that
$(\varepsilon_n)_{n\geq 0}$ is a stationary martingale difference sequence  with $\mathbf{E}[\varepsilon_i^2 | \varepsilon_1,...,\varepsilon_{i-1}]= \sigma^2$ a.s.\ for a positive constant  $\sigma.$
%%An important example is the case where  $(\varepsilon_n)_{n\geq 0}$  is a sequence of
%%i.i.d.\ random variables with mean $0 $  and finite variance $\sigma^2=\mathbf{E}\varepsilon_1^2.$
We can estimate the unknown parameter $\theta$ by the least-squares estimator given by
$$\widehat{\theta}_n= \frac{\sum_{i=1}^n Y_i Y_{i+1}}{ \sum_{i=1}^n Y_i^2}.  $$
 It is well known that $(\widehat{\theta}_n-\theta) \sqrt{ \Sigma_{i=1}^n Y_i^2 }  $ converges
 in distribution to a normal law, see Theorem 3 of Lai and Wei \cite{LW82}.
By Theorem \ref{th4}, we have the following Berry-Esseen bound for  the least-squares estimator $ \widehat{\theta}_n.$
\begin{theorem}\label{th415}
Suppose that $\mathbf{E}|\varepsilon_1|^{ 2p} < \infty$ for some  $p >1.$  If $|\theta|  <1,$ then
\begin{eqnarray}
&&\sup_{ x \in \mathbf{R}}\Big|\mathbf{P}\Big( (\widehat{\theta}_n-\theta) \sqrt{ \Sigma_{i=1}^n Y_i^2 } \leq x \sigma \Big)-\Phi \left( x\right) \Big| \nonumber \\
&&\ \ \ \ \ \ \ \ \ \ \ \ \ \ \   =  O\Bigg( n^{1- p } +  n^{ - p } \mathbf{E}\bigg|\sum_{i=1}^n   (Y_i ^{ 2}-\mathbf{E} Y_i^2)   \bigg|^{p}     \Bigg)^{1/( 2p+1) } ,
\end{eqnarray}
where $Y_n= \sum_{i=1}^n   \theta^{n-i}\varepsilon_{i } .$
\end{theorem}

\section{Proofs of theorems}\label{sec3}
\subsection{Proof  of Theorem   \ref{th1} }
We assume that   $N_n \leq 1$. Otherwise, $(\ref{ineq10})$ is trivial.

 %%since $N_n >0, $  we take $C_p$ large enough such that $C_p \, N_n^{1/(2p+1) }\geq 1.$ Then the theorem holds obviously.

Firstly, we give a lower bound for $\mathbf{P}(S_n \leq x \sqrt{[S]_n})-  \Phi(x ) , x\leq0.$ Let $\varepsilon_n \in (0, 1/2]$ be a positive number, whose exact value will be chosen later. It is easy to see that for $x\leq0,$
\begin{eqnarray}
\mathbf{P}(S_n \leq x \sqrt{[S]_n})-  \Phi(x )   &\geq & \mathbf{P}(S_n \leq x \sqrt{[S]_n}, \ [S]_n < 1+ \varepsilon_n)-  \Phi(x )  \nonumber  \\
&\geq& \mathbf{P}(S_n \leq x \sqrt{1+ \varepsilon_n }, \ [S]_n < 1+ \varepsilon_n)-   \Phi(x )\nonumber \\
&\geq& \mathbf{P}(S_n \leq x \sqrt{1+ \varepsilon_n } )- \mathbf{P}( [S]_n \geq  1+ \varepsilon_n) - \Phi(x )  \nonumber \\
& =&I_1 + I_2 - I_3, \label{fgh364}
\end{eqnarray}
where
\begin{eqnarray}
I_1&=& \mathbf{P}(S_n \leq x \sqrt{1+ \varepsilon_n } ) -  \Phi(x \sqrt{1+ \varepsilon_n }),   \nonumber  \\
I_2 &=&  \Phi(x \sqrt{1+ \varepsilon_n }) -  \Phi(x ),  \nonumber \\
I_3 &=& \mathbf{P}( [S]_n \geq  1+ \varepsilon_n). \nonumber
\end{eqnarray}
Next, we estimate $I_1, I_2$ and $I_3.$   By Haeusler's inequality \cite{H88} (see also (\ref{sfdf}) when $p\in (1, 2]$), we get the following estimation for $I_1:$
\begin{eqnarray}
I_1 \geq -C_{p,1}\Big( \sum_{i=1}^n \mathbf{E}  |X_i|^{2p}   +\mathbf{E} \big| \langle S\rangle_n-1\big|^p  \Big)^{1/(2p+1) }. \label{21b3w}
\end{eqnarray}
By one-term Taylor's expansion, we have the following estimation for $I_2:$
\begin{eqnarray}
I_2  &\geq & -c_1  e^{-x^2/2} |x| (\sqrt{1+ \varepsilon_n } -1 \big)   \nonumber \\
& \geq & -c_2 \,  \varepsilon_n. \label{21bsf3w}
\end{eqnarray}
For $I_3,$  by Markov's inequality, it follows that
\begin{eqnarray}
I_3  &=& \mathbf{P}([S]_n -\langle S\rangle_n+\langle S\rangle_n -1\geq   \varepsilon_n )   \nonumber  \\
&\leq & \mathbf{P}\Big([S]_n -\langle S\rangle_n \geq   \frac{ \varepsilon_n}2 \Big ) +\mathbf{P}\Big(\langle S\rangle_n -1\geq  \frac{ \varepsilon_n}2  \Big)   \nonumber \\
& \leq & c_{3 }\varepsilon_n^{-p} \Big(  \mathbf{E} |[S]_n -\langle S\rangle_n|^p   +    \mathbf{E} |\langle S\rangle_n -1|^p   \Big) .\label{fg563s15}
\end{eqnarray}
We distinguish two cases to estimate $I_3$.  Notice that $([S]_i - \langle S\rangle_i, \mathcal{F}_{i})_{i=0,...,n} $ is also a martingale.

\emph{Case 1}: If $p \in (1, 2],$ by the inequality  of von Bahr-Esseen \cite{V65}, it follows that
\begin{eqnarray}
 \mathbf{E} |[S]_n -\langle S\rangle_n|^p  &\leq& c_4  \sum_{i=1}^n\mathbf{E} |X_i^2   -\mathbf{E}[X_i^2 |\mathcal{F}_{i-1}] |^p  \nonumber\\
 &\leq & 2c_4  \sum_{i=1}^n \mathbf{E}[ |X_i  |^{2p}+|\mathbf{E}[X_i^2 |\mathcal{F}_{i-1}] |^p\, ] \nonumber \\
 &\leq & c_5  \sum_{i=1}^n\mathbf{E} |X_i  |^{2p} .\label{21315w}
\end{eqnarray}
Returning to (\ref{fg563s15}), we have
\begin{eqnarray}\label{213gswgd}
I_3   \leq   c_{6 }\ \varepsilon_n^{-p}\Big(  \sum_{i=1}^n\mathbf{E} |X_i  |^{2p}   +  \mathbf{E} |\langle S\rangle_n -1|^p   \Big) .
\end{eqnarray}

\emph{Case 2}: If $p >2,$ by  Rosenthal's inequality (cf.\ Theorem 2.12 of Hall  and Heyde \cite{HH80}),
 we have
\begin{eqnarray}
 \mathbf{E} | [S]_n - \langle S\rangle_n |^p    \leq  C_{p,2} \Big(  \mathbf{E}\Big( \sum_{i=1}^n \mathbf{E}[X_i^4 |\mathcal{F}_{i-1}  ]  \Big)^{p/2}  +  \sum_{i=1}^n \mathbf{E} |X_i|^{2p}  \Big ).\label{213w}
\end{eqnarray}
Noting that $X_i^4 = (X_i^2)^{(p-2)/(p-1)} (|X_i|^{2p})^{1/(p-1)}$ for $p>2$, we have by
H\"{o}lder's inequality
$$ \mathbf{E}[X_i^4|\mathcal{F}_{i-1}] \leq  \Big( \mathbf{E}[|X_i|^{2p}|\mathcal{F}_{i-1}] \Big)^{1/(p-1)}\Big(  \mathbf{E}[X_i^{2}|\mathcal{F}_{i-1}]  \Big)^{(p-2)/(p-1) } ,$$
and hence
\begin{eqnarray}
 \sum_{i=1}^n\mathbf{E}[X_i^4|\mathcal{F}_{i-1}] &\leq& \sum_{i=1}^n \Big( \mathbf{E}[|X_i|^{2p}|\mathcal{F}_{i-1}] \Big)^{1/(p-1)}\Big(  \mathbf{E}[X_i^{2}|\mathcal{F}_{i-1}]  \Big)^{(p-2)/(p-1) } \nonumber \\
&\leq&\Big(\sum_{i=1}^n \mathbf{E}[|X_i|^{2p}|\mathcal{F}_{i-1}] \Big)^{1/(p-1)}\Big(\langle S\rangle_n \Big)^{(p-2)/(p-1) }. \nonumber
\end{eqnarray}
By the inequality
 $$(a +b)^q \leq 2^q(a^q +b^q),\ \ \ \ \ a, b \geq 0 \ \ \textrm{and} \ \  q >0 ,$$
 and the fact that $(p^2-2p)/(2p-2) \leq p,$
it follows that for $p >2,$
\begin{eqnarray}
  \Big(\sum_{i=1}^n\mathbf{E}[X_i^4|\mathcal{F}_{i-1}] \Big)^{p/2}
&\leq& \Big(\sum_{i=1}^n \mathbf{E}[|X_i|^{2p}|\mathcal{F}_{i-1}] \Big)^{p/(2p-2)}\Big( \langle S\rangle_n \Big)^{(p^2-2p)/(2p-2) }  \label{ineq34} \\
&\leq& 2^p \Big(\sum_{i=1}^n \mathbf{E}[|X_i|^{2p}|\mathcal{F}_{i-1}] \Big)^{p/(2p-2)}\Big( 1+ | \langle S \rangle_n-1| ^{(p^2-2p)/(2p-2) }   \Big)   \nonumber \\
&\leq& 2^p    \Big(\sum_{i=1}^n \mathbf{E}[|X_i|^{2p}|\mathcal{F}_{i-1}] \Big)^{p/(2p-2)} \nonumber \\
 && \,  + \, 2^p    \Big(\sum_{i=1}^n \mathbf{E}[|X_i|^{2p}|\mathcal{F}_{i-1}] \Big)^{p/(2p-2)}     | \langle S \rangle_n-1| ^{p(p -2 )/(2p-2) }. \nonumber
\end{eqnarray}
As to the second term on the r.h.s.\ of the last inequality, we  use the inequality $$x^ay^{1-a} \leq x +y, \ \ \ \ \ x, y \geq0 \ \textrm{and}\  a \in [0, 1],$$
and hence
\begin{eqnarray}
  \Big(\sum_{i=1}^n\mathbf{E}[X_i^4|\mathcal{F}_{i-1}] \Big)^{p/2}  &\leq& 2^p    \Big(\sum_{i=1}^n \mathbf{E}[|X_i|^{2p}|\mathcal{F}_{i-1}] \Big)^{p/(2p-2)} \nonumber \\
 && \,   + \, 2^p    \Big(\sum_{i=1}^n \mathbf{E}[|X_i|^{2p}|\mathcal{F}_{i-1}] +   | \langle S \rangle_n-1|^p  \Big). \nonumber
\end{eqnarray}
Thus,
\begin{eqnarray}
&& \mathbf{E} \Big(\sum_{i=1}^n\mathbf{E}[X_i^4|\mathcal{F}_{i-1}] \Big)^{p/2}  \nonumber \\
&&\leq 2^p  \Bigg[ \mathbf{E}  \Big(\sum_{i=1}^n \mathbf{E}[|X_i|^{2p}|\mathcal{F}_{i-1}] \Big)^{p/(2p-2)}   +  \sum_{i=1}^n \mathbf{E} |X_i|^{2p}   +  \mathbf{E}   | \langle S \rangle_n-1|^p  \Bigg]  \nonumber \\
&&\leq 2^p  \Bigg[   \Big(\sum_{i=1}^n \mathbf{E} |X_i|^{2p}   \Big)^{p/(2p-2)}   +  \sum_{i=1}^n \mathbf{E} |X_i|^{2p}   +  \mathbf{E}   | \langle S \rangle_n-1|^p   \Bigg]. \label{fs2nl}
\end{eqnarray}
Returning to (\ref{213w}), we get for $p >2,$
\begin{eqnarray*}
 \mathbf{E} | [S]_n - \langle S\rangle_n |^p  \leq  C_{p,3} \Bigg(  \Big(\sum_{i=1}^n \mathbf{E} |X_i|^{2p}   \Big)^{p/(2p-2)}   +  \mathbf{E}   | \langle S \rangle_n-1|^p  \Bigg ).
\end{eqnarray*}
From (\ref{fg563s15}) and the last inequality,  we obtain  for $p >2,$
\begin{eqnarray}\label{213sfdw}
I_3   \leq   C_{p,4}\ \varepsilon_n^{-p} \Bigg(  \Big(\sum_{i=1}^n \mathbf{E} |X_i|^{2p}   \Big)^{p/(2p-2)}   +  \mathbf{E}   | \langle S \rangle_n-1|^p  \Bigg ).
\end{eqnarray}

By the inequalities (\ref{213gswgd}) and (\ref{213sfdw}), we always have  for $p>1,$
\begin{eqnarray}\label{213gsw}
I_3   \leq   C_{p, 5 }\ \varepsilon_n^{-p} \Bigg( \sum_{i=1}^n \mathbf{E} |X_i|^{2p}  + \Big(\sum_{i=1}^n \mathbf{E} |X_i|^{2p}   \Big)^{p/(2p-2)}   +  \mathbf{E}   | \langle S \rangle_n-1|^p  \Bigg ).
\end{eqnarray}
Combining (\ref{fgh364}), (\ref{21b3w}), (\ref{21bsf3w}) and (\ref{213gsw}) together, we deduce that  for $p >1,$
\begin{eqnarray}
&& \mathbf{P}(S_n \leq x \sqrt{[S]_n})-  \Phi(x ) \nonumber \\
 & &\geq- C_{p, 1 }\, \Big(\sum_{i=1}^n \mathbf{E}  |X_i|^{2p}  +\mathbf{E}  \big| \langle S\rangle_n-1\big|^p    \Big)^{1/(2p+1) } -\, c_2 \,\varepsilon_n\nonumber  \\
&&\  \ \ -\,  C_{p, 5 }\, \varepsilon_n^{-p}\Bigg(  \sum_{i=1}^n\mathbf{E} |X_i  |^{2p}  + \Big(\sum_{i=1}^n \mathbf{E} |X_i|^{2p}   \Big)^{p/(2p-2)} +  \mathbf{E} |\langle S\rangle_n -1|^p   \Bigg) . \nonumber
\end{eqnarray}
Taking
\begin{eqnarray}\label{fhk02l}
\varepsilon_n=\Bigg(  \sum_{i=1}^n\mathbf{E} |X_i  |^{2p}   +\Big(\sum_{i=1}^n \mathbf{E} |X_i|^{2p}   \Big)^{p/(2p-2)}+  \mathbf{E} |\langle S\rangle_n -1|^p   \Bigg)^{ 1/(p+1)  },
\end{eqnarray}
 we obtain  for $x\leq 0$ and  $p >1,$
\begin{eqnarray}
&&\mathbf{P}(S_n \leq x \sqrt{[S]_n})-  \Phi(x ) \nonumber \\
&&\geq -C_{p, 1 }\Big(\sum_{i=1}^n \mathbf{E}  |X_i|^{2p}  +\mathbf{E}  \big| \langle S\rangle_n-1\big|^p    \Big)^{1/(2p+1) } \nonumber  \\
&& \ \ \   - \,C_{p, 6}\,  \Bigg(  \sum_{i=1}^n\mathbf{E} |X_i  |^{2p}   +\Big(\sum_{i=1}^n \mathbf{E} |X_i|^{2p}   \Big)^{p/(2p-2)} +  \mathbf{E} |\langle S\rangle_n -1|^p   \Bigg)^{1/(p+1) }  \nonumber  \\
&&\geq -C_{p,7}\Big(\sum_{i=1}^n \mathbf{E}  |X_i|^{2p}  +\mathbf{E}  \big| \langle S\rangle_n-1\big|^p    \Big)^{1/(2p+1) } ,\label{fgjk16}
\end{eqnarray}
where the last line follows from the fact that $p/((2p-2)(p+1)) \geq 1 /(2p+1)$ and $N_n \leq 1$.

Secondly, we give an upper bound for $\mathbf{P}(S_n \leq x \sqrt{[S]_n})- \Phi(x )  , x\leq 0.$ It is obvious that for $x\leq 0,$
\begin{eqnarray}
&& \mathbf{P}(S_n \leq x \sqrt{[S]_n})-  \Phi(x ) \nonumber \\
& &\leq \mathbf{P}(S_n \leq x \sqrt{[S]_n}, \ [S]_n > 1- \varepsilon_n)-  \Phi(x )\nonumber \\
 &&\ \ \  \, +\, \mathbf{P}(S_n \leq x \sqrt{[S]_n}, \ [S]_n \leq 1- \varepsilon_n) \nonumber \\
&&\leq\mathbf{P}(S_n \leq x \sqrt{1- \varepsilon_n }, \ [S]_n > 1-\varepsilon_n)-  \Phi(x )+   \mathbf{P}(  [S]_n \leq 1- \varepsilon_n)   \nonumber \\
&&\leq \mathbf{P}(S_n \leq x \sqrt{1- \varepsilon_n } )-\Phi(x\sqrt{1- \varepsilon_n } )+ \Phi(x \sqrt{1- \varepsilon_n })-\Phi(x ) \nonumber \\
&& \ \ \ \, + \,   \mathbf{P}(  [S]_n \leq 1- \varepsilon_n) \nonumber  \\
&&=I_4 + I_5 +I_6  . \nonumber
\end{eqnarray}
Following the same lines as in the  proof of (\ref{fgjk16}),  we get for $x\leq 0$ and  $p >1,$
\begin{eqnarray} \label{ineq24}
\mathbf{P}(S_n \leq x \sqrt{[S]_n})-  \Phi(x )  \, \leq \,   C_{p,8}\Big(\sum_{i=1}^n \mathbf{E}  |X_i|^{2p}  +\mathbf{E}  \big| \langle S\rangle_n-1\big|^p \Big)^{1/(2p+1) } .
\end{eqnarray}
Combining (\ref{fgjk16}) and (\ref{ineq24}) together, we get for $p >1,$
\begin{eqnarray}\label{ineq25}
\sup_{x\leq 0}\Big|\mathbf{P}(S_n \leq x \sqrt{[S]_n})-  \Phi(x ) \Big| \, \leq \,   C_{ p, 8}\Big(\sum_{i=1}^n \mathbf{E}  |X_i|^{2p}  +\mathbf{E}  \big| \langle S\rangle_n-1\big|^p    \Big)^{1/(2p+1) } .
\end{eqnarray}

 Notice that $(-S_k,\mathcal{F}_k)_{k=0,...,n}$ is also a martingale. Applying   the last inequality to $(-S_k,\mathcal{F}_k)_{k=0,...,n}$,
  we get
 \begin{eqnarray}
&&\sup_{x> 0}\Big|\mathbf{P}(S_n \leq  x \sqrt{[S]_n})-  \Phi( x ) \Big| \nonumber \\
  &&= \sup_{x> 0}\Big|\mathbf{P}( S_n  \leq x \sqrt{[S]_n})-1+1-  \Phi(x ) \Big| \nonumber \\
&&= \sup_{x> 0}\Big| \Phi(-x )-\mathbf{P}(- S_n <-x \sqrt{[S]_n})  \Big|
\nonumber \\
&&\leq C_{ p,9}\Big(\sum_{i=1}^n \mathbf{E}  |X_i|^{2p}  +\mathbf{E}  \big| \langle S\rangle_n-1\big|^p   \Big)^{1/(2p+1) } . \label{ineq66}
\end{eqnarray}
Combining the inequalities (\ref{ineq25}) and (\ref{ineq66}) together,  we obtain
 \begin{eqnarray}
\sup_{x\in \mathbf{R}}\Big|\mathbf{P}(S_n \leq  x \sqrt{[S]_n})-  \Phi( x ) \Big|
 \leq  C_{p,10}\Big(\sum_{i=1}^n \mathbf{E}  |X_i|^{2p}  + \mathbf{E} \big| \langle S\rangle_n-1\big|^p    \Big)^{1/(2p+1) } , \label{ineq6fs}
\end{eqnarray}
which gives the desired inequality (\ref{ineq10}).

Next we give a proof of (\ref{ineq07}). We follow the example  of Haeusler \cite{H88}.
Let $(\alpha_n)_{n\geq 1}$ be a sequence of positive numbers such that $\alpha_n \rightarrow 0$ as $n\rightarrow \infty.$
Define  function $f_n: \mathbf{R} \rightarrow [0, \, \infty)$ as follows
\begin{displaymath}
f_n(x)= \left\{ \begin{array}{ll}
x^{-1}, & \textrm{if $ \frac12\sqrt{\alpha_n}   \leq  x < \infty$,}\\
0, & \textrm{otherwise} .
\end{array} \right.
\end{displaymath}
Furthermore, let $X_1, ..., X_{n-1}$ be independent and normally distributed random variables  with mean $0$ and
variance $1/(n-1).$ Denote $\nu_x $ the one-point mass concentration at $x.$   Define the random variable $X_n$
 such that its conditional distribution, given  $X_1, ..., X_{n-1}$,  is
$$ \mathbf{P}(X_n \in  \cdot | S_{n-1}=x) = \frac12 \nu_{- \alpha_n f_n(x)}(\cdot) + \frac12 \nu_{  \alpha_n f_n(x)}(\cdot), $$
where $S_{n-1}=\sum_{i=1}^{n-1}X_i.$  Denote  $\mathcal{F}_{i}$ the natural filtration of $X_1, ..., X_{n}, $ that is $\mathcal{F}_0$
being the trivial $\sigma-$field and  $\mathcal{F}_{i}= \sigma \{ X_1, ..., X_{i} \}, i=1,...,n.$ Clearly, $(X_i, \mathcal{F}_{i})_{i=0,...,n}$ is
a finite sequence of martingale differences.  Moreover, it holds
$$  \sum_{i=1}^{n-1} \mathbf{E} |X_i|^{2p}   = \frac{n-1}{(n-1)^{p}}   \mathbf{E}  |  \mathcal{N} (0, 1) |^{2p}   \sim C_{p, 11 } n^{1-p}  $$
and
\begin{eqnarray}
   \mathbf{E}  |X_n|^{2p}   &=&    \int_{-\infty}^{\infty} \int_{-\infty}^{\infty} |  y|^{2p} \mathbf{P}( X_n \in dy |  \mathcal{N}(0 , 1 ) =x )   \mathbf{P}( \mathcal{N}(0 , 1 ) \in dx)  \nonumber \\
   & =&  \frac{1}{\sqrt{2 \pi} } \int_{\frac12 \sqrt{\alpha_n}  }^\infty \Big|  \frac{\alpha_n}{x }  \Big |^{2p}  e^{-\frac12x^2 } dx  \nonumber \\
   &\sim&C_{p, 12} \, \alpha_n^{p+ \frac1 2 }
\end{eqnarray}
for some constants $0<C_{p, 11}, C_{p, 12} < \infty,$  where $\mathcal{N}(0 , 1 )$ is a standard random variable. Similarly, we have
$$\sum_{i=1}^{n-1} \mathbf{E}[X_i^{2}| \mathcal{F}_{i-1} ] = \sum_{i=1}^{n-1} \mathbf{E} X_i^{2}    = 1 $$
and
\begin{eqnarray}
  \mathbf{E}  \big| \langle S\rangle_n-1\big|^p   &=&\mathbf{E}  \big| \mathbf{E}[X_n^{2}| \mathcal{F}_{n-1} ]\big|^p   \nonumber \\
  & =&    \int_{-\infty}^{\infty}  \Big|  \int_{-\infty}^{\infty} y^2 \mathbf{P}( X_n \in dy |  \mathcal{N}(0 , 1 ) =x ) \Big|^{p}   \mathbf{P}( \mathcal{N}(0 , 1 ) \in dx)\nonumber \\
  &=& \frac{1}{\sqrt{2 \pi} } \int_{ \frac12 \sqrt{\alpha_n}  }^\infty \Big|  \frac{\alpha_n}{x }  \Big |^{2p}  e^{-\frac12 x^2 } dx   \nonumber \\
  &\sim& C_{p,12} \, \alpha_n^{p+ \frac1 2 }.
\end{eqnarray}
Thus
$$N_n \sim   C_p \, \alpha_n^{p+ \frac1 2 } ,$$
 where $C_p$ is a positive constant  depending only on $p$.
On the other hand, we  have
\begin{eqnarray}
  \mathbf{P}\Big( \frac{S_n }{\sqrt{[S]_n}} \leq 0 \Big)
  & =&   \mathbf{P}\Big( X_n + S_{n-1}\leq  0  \Big)  \nonumber \\
    &=& \int_{-\infty}^{\infty}\mathbf{P}( X_n \leq  -x |  \mathcal{N}(0 , 1 ) =x )       \mathbf{P}( \mathcal{N}(0 , 1 ) \in dx)\nonumber \\
    &=& \int_{-\infty}^{\frac12 \sqrt{\alpha_n} }\mathbf{P}( X_n \leq  -x |  \mathcal{N}(0 , 1 ) =x )       \mathbf{P}( \mathcal{N}(0 , 1 ) \in dx)\nonumber \\
    &&  + \int_{\frac12 \sqrt{\alpha_n}}^{\infty}\mathbf{P}( X_n \leq  -x |  \mathcal{N}(0 , 1 ) =x )       \mathbf{P}( \mathcal{N}(0 , 1 ) \in dx)\nonumber \\
  &=& \frac{1}{\sqrt{2 \pi} } \int_{ -\infty }^0   e^{-\frac12 x^2 } dx  + \frac{1}{\sqrt{2 \pi} } \int_{\frac12 \sqrt{\alpha_n}}^{\sqrt{\alpha_n}} \frac12 \,  e^{-\frac12x^2 } dx   \nonumber \\
  &=&\Phi(0) +    \frac{1}{ 4 \sqrt{2 \pi} }   \sqrt{\alpha_n}  \Big(1+o(1) \Big).  \nonumber
\end{eqnarray}
Hence, we deduce that
\begin{eqnarray}
 &&\sup_{ x \in \mathbf{R}}\Big|\mathbf{P}\Big( \frac{S_n}{\sqrt{[S]_n}}  \leq x \Big)-\Phi \left( x\right) \Big|  \, N_n^{-1/(2p+1) } \nonumber \\
  && \geq \Big|\mathbf{P}\Big( \frac{S_n}{\sqrt{[S]_n}}  \leq 0\Big)-\Phi \left( 0\right) \Big|  \, (  C_p \, \alpha_n^{p+ \frac12  } )^{-1/(2p+1) }  \Big(1+o(1) \Big)   \nonumber\\
  &&=\frac{1}{ 4 \sqrt{2 \pi} }   \sqrt{\alpha_n}\, (  C_p \, \alpha_n^{p+ \frac 12  } )^{-1/(2p+1) }  \Big(1+o(1) \Big)    \nonumber\\
  && \sim \frac{1}{ 4 \sqrt{2 \pi} }  \, (  C_p \,  )^{-1/(2p+1) } .    \nonumber
\end{eqnarray}
This completes the proof of Theorem \ref{th1}.

\subsection{Proof  of Theorem   \ref{th4} }
First, we give a lower bound for $\mathbf{P}(S_n \leq x \sqrt{\langle S\rangle_n})-  \Phi(x ) , x\leq0.$ Let $\varepsilon_n \in (0, 1/2]$ be a positive number, whose exact value will be chosen later. It is easy to see that for $x\leq0,$
\begin{eqnarray}
\mathbf{P}\big(S_n \leq x \sqrt{\langle S\rangle _n} \big)-  \Phi(x )   &\geq & \mathbf{P}(S_n \leq x \sqrt{\langle S\rangle _n}, \ \langle S\rangle _n< 1+ \varepsilon_n)-  \Phi(x )  \nonumber  \\
&\geq& \mathbf{P}(S_n \leq x \sqrt{1+ \varepsilon_n }, \ \langle S\rangle _n < 1+ \varepsilon_n)-   \Phi(x )\nonumber \\
&\geq& \mathbf{P}(S_n \leq x \sqrt{1+ \varepsilon_n } )- \mathbf{P}( \langle S\rangle _n \geq  1+ \varepsilon_n) - \Phi(x )  \nonumber \\
& =& I_1 + I_2 - I_7,  \nonumber
\end{eqnarray}
where
\begin{eqnarray}
I_7  =  \mathbf{P}( \langle S\rangle_n \geq  1+ \varepsilon_n). \nonumber
\end{eqnarray}
For $I_7,$  by Markov's inequality, it follows that
\begin{eqnarray}\label{21qf3w}
I_7   \leq    \varepsilon_n^{-p} \mathbf{E} |\langle S\rangle_n -1|^p  .
\end{eqnarray}
Combining the estimations (\ref{21b3w}), (\ref{21bsf3w}) and (\ref{21qf3w}) together,  we have for $p >1,$
\begin{eqnarray}
&& \mathbf{P}(S_n \leq x \sqrt{[S]_n})-  \Phi(x ) \nonumber \\
 & &\geq- C_{p, 1 }\Big(\sum_{i=1}^n \mathbf{E}  |X_i|^{2p}  +\mathbf{E} \big| \langle S\rangle_n-1\big|^p    \Big)^{1/(2p+1) } -\, c_2 \,\varepsilon_n  \ -\,   \varepsilon_n^{-p}\mathbf{E} |\langle S\rangle_n -1|^p  . \nonumber
\end{eqnarray}
Next we carry out an argument as the proof of Theorem \ref{th1} with
\begin{eqnarray}
\varepsilon_n=\Big(   \mathbf{E} |\langle S\rangle_n -1|^p   \Big)^{ 1/(p+1)  },
\end{eqnarray}
what we obtain is
\begin{eqnarray}
 \sup_{ x \in \mathbf{R}}\Big|\mathbf{P}\Big( \frac{S_n}{\sqrt{\langle S \rangle_n}}  \leq x \Big)-\Phi \left( x\right) \Big|
 \leq  C_{p,2}\Big(\sum_{i=1}^n \mathbf{E}  |X_i|^{2p}  +\mathbf{E}  \big| \langle S\rangle_n-1\big|^p   \Big)^{1/(2p+1) } ,
\end{eqnarray}
that is inequality (\ref{ineq10})   holds when
 $ S_n/\sqrt{[S]_n} $ is replaced by $ S_n/\sqrt{\langle S \rangle_n} .$
The proof of optimality is similar to the proof of (\ref{ineq07}).
This completes the proof of Theorem \ref{th4}.

\subsection{Proof  of Theorem \ref{th415} }
The proof of theorem is based on Theorem \ref{th4}.
It is easy to see that
$$ \frac1{\sigma}\,(\widehat{\theta}_n-\theta) \sqrt{  \Sigma_{i=1}^n Y_i^2 }= \frac{\sum_{i=1}^n Y_i \varepsilon_{i+1}}{ \sigma \sqrt{ \sum_{i=1}^n Y_i^2}}  .    $$
Notice that $Y_n= \sum_{i=1}^n   \theta^{n-i}\varepsilon_{i } .$
Set
$$X_i= \frac{  Y_i \varepsilon_{i+1}}{ \sigma \sqrt{ \sum_{i=1}^n \mathbf{E}Y_i^2 }}\ \ \ \ \ \textrm{and} \ \ \ \  \mathcal{F}_i=\sigma\{ \varepsilon_k,\  1\leq k \leq i+1  \}. $$
Then it is easy to see that $(X_i, \mathcal{F}_i)_{i=0,...,n}$ is a sequence of martingale differences,
 and that
 $$\frac1{\sigma}\,(\widehat{\theta}_n-\theta) \sqrt{  \Sigma_{i=1}^n Y_i^2 }= \frac{S_n}{ \sqrt{\langle S\rangle_n}  }.$$
Moreover, we have
$$\mathbf{E}Y_n^2  = \sum_{i=1}^n   \theta^{2(n-i)} \sigma^2= \frac{1- \theta^{2\, n} }{1-\theta^2 }  \sigma^2 $$ % \sim  \frac{1 }{1-\theta^2 }  \sigma^2 $$
and
$$\sum_{i=1}^n \mathbf{E}Y_i^2 = \sum_{i=1}^n \frac{1- \theta^{2\, i} }{1-\theta^2 }  \sigma^2.  $$
By  Rosenthal's inequality, we also have
\begin{eqnarray}
  \mathbf{E}|Y_n|^{2p}  &\leq& C_p \Big( (  \mathbf{E}Y_n^2)^{p} + \sum_{i=1}^n \mathbf{E}|\theta^{n-i}\varepsilon_{i }|^{2p} \Big) \nonumber \\
&\leq& C_p \Bigg(  \Big(  \frac{1- \theta^{2\, n} }{1-\theta^2 }  \Big)^{p } \sigma^{2p} + \frac{1- |\theta|^{2p\, n} }{1-|\theta|^{2p} } \mathbf{E} |\varepsilon_{1 }|^{2p} \Bigg) \nonumber %\\   &\leq&  C_p \Big( \frac{1  }{ ( 1-\theta^2)^{p/2}  }  \sigma^p + \frac{1  }{1-|\theta|^p } \mathbf{E} |\varepsilon_{1 }|^p \Big)
\end{eqnarray}
and
\begin{eqnarray*}
  \sum_{i=1}^n \mathbf{E}|Y_i|^{2p }  \leq     C_p  \, \sum_{i=1}^n \Bigg(  \Big(  \frac{1- \theta^{2\, i} }{1-\theta^2 }  \Big)^{p } \sigma^{2p} + \frac{1- |\theta|^{2p\, i} }{1-|\theta|^{2p} } \mathbf{E} |\varepsilon_{1 }|^{2p } \Bigg).
\end{eqnarray*}
Thus
\begin{eqnarray}
 \frac{\sum_{i=1}^n \mathbf{E} |X_i|^{2p} }{ (\sum_{i=1}^n \mathbf{E} X_i^2 )^{p }} \leq  C_p  \, \sum_{i=1}^n \Bigg(  \Big(  \frac{1- \theta^{2\, i} }{1-\theta^2 }  \Big)^{p } \sigma^{2p } + \frac{1- |\theta|^{2p\, i} }{1-|\theta|^{2p} } \mathbf{E} |\varepsilon_{1 }|^{2p} \Bigg) \Big / \Big(\sum_{i=1}^n \frac{1- \theta^{2\, i} }{1-\theta^2 }  \sigma^2 \Big)^{p}. \label{fsf}
\end{eqnarray}
%For $\theta,$ we distinguish two cases.
If $|\theta| < 1,$ inequality (\ref{fsf}) implies that
\begin{eqnarray}
 \frac{\sum_{i=1}^n \mathbf{E} |X_i|^{ 2p} }{ (\sum_{i=1}^n \mathbf{E} X_i^2 )^{p }} &\leq&  C_p   \, \sum_{i=1}^n \Bigg(  \frac{   \sigma^{2p} }{\big(1-\theta^2\big)^{p } } + \frac{  \mathbf{E} |\varepsilon_{1 }|^{2p} }{1-|\theta|^{2p} } \Bigg) \Big / \Big(\sum_{i=1}^n    \sigma^2 \Big)^{p } \nonumber\\
 & \leq & C_p \Bigg(  \frac{ 1 }{\big(1-\theta^2\big)^{p } } + \frac{1 }{1-|\theta|^{2p} }  \frac{ \mathbf{E} |\varepsilon_{1 }|^{2p}}{ \sigma^{2p}}  \Bigg)  n^{1- p }.\nonumber
\end{eqnarray}
It is obvious that
\begin{eqnarray}
 \langle S\rangle_n =\frac{\sum_{i=1}^n   Y_i ^{ 2} }{  \sum_{i=1}^n \mathbf{E} Y_i^2 }.
\end{eqnarray}
By Theorem \ref{th4}, we obtain
\begin{eqnarray}
&&\sup_{ x \in \mathbf{R}}\Big|\mathbf{P}\Big( (\widehat{\theta}_n-\theta) \sqrt{ \Sigma_{i=1}^n Y_i^2 } \leq x \sigma \Big)-\Phi \left( x\right) \Big| \nonumber \\
 &&\leq  C_{p, \theta} \, \Bigg( \frac{\sum_{i=1}^n \mathbf{E} |X_i|^{2p}  }{ (\sum_{i=1}^n \mathbf{E}  X_i^2  )^{p }}  + \mathbf{E}\Bigg|\frac{\sum_{i=1}^n   Y_i ^{ 2} }{  \sum_{i=1}^n \mathbf{E} Y_i^2 }-1 \Bigg|^{p }    \Bigg)^{1/( 2p+1) } \nonumber \\
&&= O\Bigg( n^{1- p } +  n^{ - p } \mathbf{E}\bigg|\sum_{i=1}^n   (Y_i ^{ 2}-\mathbf{E} Y_i^2)   \bigg|^{p  }   \Bigg)^{1/(2 p+1) } . \nonumber
\end{eqnarray}
This completes the proof of theorem.
%Thus
%\begin{eqnarray}
% \frac{\sum_{i=1}^n \mathbf{E} |X_i|^{ p} }{ (\sum_{i=1}^n \mathbf{E} X_i^2 )^{p/2}} \leq  C_p  \, \sum_{i=1}^n \Bigg(  \Big(  \frac{1- \theta^{2\, i} }{1-\theta^2 }  \Big)^{p/2} \sigma^p + \frac{1- |\theta|^{p\, i} }{1-|\theta|^p } \mathbf{E} |\varepsilon_{1 }|^p \Bigg) \Big / \Big(\sum_{i=1}^n \frac{1- \theta^{2\, i} }{1-\theta^2 }  \sigma^2 \Big)^{p/2}.
%\end{eqnarray}
%\emph{Case 2}: If $|\theta| >1,$ inequality (\ref{fsf}) implies that
%\begin{eqnarray}
% \frac{\sum_{i=1}^n \mathbf{E} |X_i|^{ p} }{ (\sum_{i=1}^n \mathbf{E} X_i^2 )^{p/2}} &\leq&  C_p   \, \sum_{i=1}^n \Bigg(  \frac{ |\theta|^{pi}  \sigma^p }{\big(\theta^2 -1\big)^{p/2} } + \frac{ |\theta|^{pi} \mathbf{E} |\varepsilon_{1 }|^p }{|\theta|^p-1} \Bigg) \Big / \Big(\sum_{i=1}^n  \theta^{2 (i-1)}  \sigma^2 \Big)^{p/2} \nonumber\\
% & \leq & C_{2,p} \Bigg(  \frac{ 1 }{\big(\theta^2 -1\big)^{p/2} } + \frac{1 }{|\theta|^p -1 }  \frac{ \mathbf{E} |\varepsilon_{1 }|^p}{ \sigma^p}  \Bigg)  n^{1- p/2}.\nonumber
%\end{eqnarray}

\subsection*{Acknowledgements}
The research is partially supported by Hong Kong RGC GRF 14302515.  Fan has been partially  supported by the
 National Natural Science Foundation of China (Grant nos.\ 11601375 and  11626250).

\section*{Referees}

\end{document}